\documentclass[sn-mathphys,Numbered]{sn-jnl}

\usepackage{graphicx}%
\usepackage{multirow}%
\usepackage{amsmath,amssymb,amsfonts}%
\usepackage{amsthm}%
\usepackage{mathrsfs}%
\usepackage[title]{appendix}%
\usepackage{xcolor}%
\usepackage{textcomp}%
\usepackage{manyfoot}%
\usepackage{booktabs}%
\usepackage{algorithm}%
\usepackage{algorithmicx}%
\usepackage{algpseudocode}%
\usepackage{listings}%

\theoremstyle{thmstyleone}%
\newtheorem{theorem}{Theorem}
\theoremstyle{thmstyletwo}%
\newtheorem{thm}{Theorem}[section]

\theoremstyle{definition}

\theoremstyle{remark}

\newtheorem{corollary}[theorem]{Corollary}
\newtheorem{lemma}[theorem]{Lemma}

\newtheorem{example}{Example}%
\newtheorem{remark}{Remark}
\theoremstyle{thmstylethree}%

\begin{document}

\title[Jensen--Mercer-type Inequalities and Strongly $f$-Divergences]{Jensen--Mercer-type Inequalities and Strongly $f$-Divergences}


\author[1]{\fnm{Vidyanshu Kumar} \sur{Rai}}\email{thakurrk@nitj.ac.in; vidyanshukr.ma.25@nitj.ac.in}
\equalcont{These authors contributed equally to this work.}

\author*[2]{\fnm{Ravinder} \sur{Kumar}}\email{thakurrk@nitj.ac.in}
\equalcont{These authors contributed equally to this work.}


\affil[1,2]{\orgdiv{Department of Mathematics and Computing}, \orgname{Dr BR Ambedkar National Institute of
		Technology Jalandhar}, \orgaddress{\street{GT Road}, \city{Jalandhar}, \postcode{144008}, \state{Punjab}, \country{India}}}

%


\abstract{In this paper, we develop refined Jensen--Mercer-type inequalities for strongly convex functions, including reverse inequalities in the differentiable case. We identify the relationships between these inequalities and known bounds. Our method exploits a structural decomposition of strongly convex functions together with classical convexity arguments. As applications, we obtain explicit and computable bounds for strongly $f$-divergences, leading to new quantitative estimates for major divergence measures such as the Kullback--Leibler divergence, $\chi^{2}$-divergence, Hellinger divergence, Bhattacharyya-type divergence, and Jeffreys divergence.}


\keywords{Strongly convex function; Jensen-Mercer inequality; Hermite-Hadamard integral inequality; Strongly $f$-divergences; Kullback-Leibler divergence}


\pacs[Mathematics Subject Classification]{Primary 39B62; Secondary 26D10}

\maketitle
\section{Introduction}\label{In}
Convex functions occupy a central position in analysis, optimization, 
probability, and information theory. Numerous refinements and 
structural extensions of convexity have been developed; 
see \cite{niculescu06} for a comprehensive account.\\Strong convexity, introduced by Polyak \cite{polyak66}, 
has proved particularly influential in optimization theory 
and in the quantitative study of inequalities. 
For recent developments and applications, we refer to 
\cite{ivelic19, ivelic24, moradi18}.\\
Let $I$ be a real interval. A function $f:I\to\mathbb{R}$ is said to be 
\emph{strongly convex} with modulus $c>0$ if for all $x,y\in I$ and 
$\lambda\in[0,1]$,
\begin{equation}\label{1}	f(\lambda x+(1-\lambda)y)\le \lambda f(x)+(1-\lambda)f(y)-c\lambda(1-\lambda)(x-y)^2.
\end{equation}
If $-f$ is strongly convex, then $f$ is said to be \emph{strongly concave}.
Since $c\lambda(1-\lambda)(x-y)^2\ge 0$, every strongly convex function is convex; however, the converse does not hold. For example, every linear function is convex but not strongly convex.\\
A well-known characterization of strongly convex functions states that if $f''(x)\ge 2c>0$, then $f$ is strongly convex with modulus $c$. Another standard characterization asserts that the function $g(x)=f(x)-cx^2$ (where $g:I\to\mathbb{R}$) is convex if and only if $f$ is strongly convex with modulus $c$; see \cite{roberts73, nikodem14, ivelic24}.\\
It was shown in \cite{Ivelic2024a} that the following lemma is an immediate consequence of the above characterization of strongly convex functions together with an application of \cite[Theorem~1.6]{pecaric92} to the convex function $g=f-c\,\mathrm{id}^2$, where $\mathrm{id}$ is the identity function.
\vspace{0.2cm}
\begin{lemma}\label{l1} Let $f:I\to\mathbb{R}$ be a differentiable strongly convex function with modulus $c>0$. Then
	\begin{equation*}
		f(x)\ge f(y)+f'(y)(x-y)+c(x-y)^2,
	\end{equation*}
	for all $x\in I$ and $y\in I^\circ$ (the interior of $I$).
\end{lemma}
A strengthened version of Jensen's inequality for strongly convex 
functions was proved in \cite{nikodem14}:
\begin{equation}\label{3}f(\bar{x})	\le	\sum_{j=1}^n a_j f(x_j)	- c\sum_{j=1}^n a_j (x_j-\bar{x})^2,
\end{equation}
where $a_j>0$, $\sum_{j=1}^n a_j=1$, $x_j\in[n,N]$, and 
$\bar{x}=\sum_{j=1}^n a_j x_j$.
In the convex case $c=0$, inequality \eqref{3} reduces to classical Jensen's inequality.
Mercer \cite{mercer03} established the Jensen--Mercer inequality
\begin{equation}\label{5}f(N+n-\bar{x})	\le	f(n)+f(N)-\sum_{j=1}^n a_j f(x_j),
\end{equation}
for convex functions $f$ on $[n,N]$. This follows from the pointwise inequality
\begin{equation}\label{6}f\left(N+n-x\right)\le f(n)+f(N)-f(x), \quad n\le x\le N.\end{equation} 
Several refinements and extensions were obtained in 
\cite{khan20,butt22,horvath21,ivelic11,khan22,krnic12,nikodem14,sayyari23,roberts73,jarad23}.\\Here, we discuss a variant of the Jensen--Mercer inequality proved in \cite{moradi18}, that is,
\begin{equation}\label{7}\begin{aligned}	f\left(N+n-\bar{x} \right) \le  f(n)+f(N)-\sum_{j=1}^{n} a_j f(x_j)\\	& \hspace{-5.5cm}- c\left[2(n-N)^2\sum_{j=1}^{n} a_j\lambda_j(1-\lambda_j)+\sum_{j=1}^{n} a_j (x_j -\bar{x})^2\right]. \end{aligned}
\end{equation}
It was also shown in \cite{moradi18} that, for a strongly convex function,
\begin{equation}\label{7*}
	f\left(N+n-x\right)\le f(N)+f(n)-f(x)-2c(N-x)(x-n), \quad n\le x\le N.
\end{equation}
The Hermite--Hadamard inequality for a strongly convex function $f:I\to\mathbb{R}$ with modulus $c>0$ is given in \cite{merentes10}:
\begin{equation}\label{8} f\!\left(\frac{n+N}{2}\right)+ \frac{c}{12}(n-N)^2\le \frac{1}{N-n}\int_n^N f(t)\,dt \le \frac{f(n)+f(N)}{2}-\frac{c}{6}(n-N)^2,
\end{equation}
where $n,N \in I$ and $n<N$. Note that 
\begin{equation*}	\begin{aligned}	\frac{1}{N-n}\int_n^N f(t)\,dt = \int_0^1 f((1-t)n+tN)\,dt = \int_0^1 f((1-t)N+tn)\,dt.	\end{aligned}
\end{equation*}
Therefore, \eqref{8} can also be written in the following form,
\begin{equation*}	\begin{aligned}	f\!\left(\frac{n+N}{2}\right)+ \frac{c}{12}(n-N)^2\le \int_0^1 f((1-t)n+tN)\,dt \le \frac{f(n)+f(N)}{2}-\frac{c}{6}(n-N)^2.	\end{aligned}
\end{equation*}
For a convex function, the corresponding inequality is  
\begin{equation}\label{he1} 	\begin{aligned}	f\!\left(\frac{n+N}{2}\right) \le \frac{1}{N-n}\int_n^N f(t)\,dt \le \frac{f(n)+f(N)}{2}.	\end{aligned}
\end{equation}
This is the classical Hermite--Hadamard inequality; see \cite{cerone10}.\\
Throughout this paper, $I$ denotes an open interval containing $[n,N]$. Unless otherwise stated, all functions are assumed to be defined on $I$.\\

The aim of this paper is to incorporate strong convexity into Jensen--Mercer and Hermite--Hadamard type inequalities and to apply the resulting refinements to strong $f$-divergences. The main novelty of the paper lies in deriving explicit quadratic correction terms involving the modulus of strong convexity. As an application, we establish new bounds for several classical divergence measures, including the Kullback--Leibler, $\chi^2$, Hellinger, Bhattacharyya, and Jeffreys divergences.\\
The paper is organized as follows. In Section~\ref{mr}, we establish refined Jensen--Mercer and Hermite--Hadamard type inequalities for strongly convex functions. Section~\ref{Ap} is devoted to applications of these results to strong $f$-divergences.
\section{Refined Jensen--Mercer Type Inequalities}\label{mr}
In this section, we establish strengthened versions of Jensen--Mercer type inequalities for strongly convex functions.

\vspace{0.2cm}
\begin{thm}\label{JMstrong} Let $f:I\to\mathbb{R}$ be a strongly convex function with modulus $c>0$. Let $n<N$, $n\le x_j\le N$, $a_j>0$ and $\sum_{j=1}^k a_j=1$. Denote $\bar{x}=\sum_{j=1}^k a_jx_j$. Then
\begin{equation}\label{JMrefined}
f(N+n-\bar{x})\le f(N)+f(n)-\sum_{j=1}^k a_jf(x_j)
-c\left[2(N-\bar{x})(\bar{x}-n)-\sum_{j=1}^k a_j(x_j-\bar{x})^2\right].
\end{equation}\end{thm}
\begin{proof}
	Since $f$ is strongly convex with modulus $c$, the function
	$g(t)=f(t)-ct^2$ is convex on $[n,N]$. Applying the classical
	Jensen--Mercer inequality to $g$, we obtain
	$$	g(N+n-\bar{x})
	\le g(N)+g(n)-\sum_{j=1}^k a_jg(x_j).	$$
	Substituting $g(t)=f(t)-ct^2$ and rearranging gives
	\begin{align*}
		f(N+n-\bar{x})
		&\le f(N)+f(n)-\sum_{j=1}^k a_jf(x_j)\\
		&\quad-c\left[
		N^2+n^2-(N+n-\bar{x})^2
		-\sum_{j=1}^k a_jx_j^2
		\right].
	\end{align*}
	Since
	\[
	\sum_{j=1}^k a_jx_j^2
	=\bar{x}^{\,2}
	+\sum_{j=1}^k a_j(x_j-\bar{x})^2,
	\]
	the expression inside the square brackets becomes
	\begin{align*}
		&N^2+n^2-(N+n-\bar{x})^2
		-\bar{x}^{\,2}
		-\sum_{j=1}^k a_j(x_j-\bar{x})^2\\
		&\quad=2\bigl((N+n)\bar{x}-Nn-\bar{x}^{\,2}\bigr)
		-\sum_{j=1}^k a_j(x_j-\bar{x})^2\\
		&\quad=2(N-\bar{x})(\bar{x}-n)
		-\sum_{j=1}^k a_j(x_j-\bar{x})^2.
	\end{align*}
	Substituting this expression into the preceding inequality proves
	\eqref{JMrefined}.
\end{proof}
We now present the following theorem, which extends the Jensen--Mercer inequality to two weighted means.\vspace{0.2cm}
\begin{thm}\label{t1}
	Let $f:I\to\mathbb{R}$ be a strongly convex function with modulus $c>0$.
	Let $n<N$, $n\le x_j,y_j\le N$ $(j=1,2,\ldots,k)$, and let $a_j>0$
	such that
	\begin{align*}
		\sum_{j=1}^k a_j=1,\qquad
		\bar{x}=\sum_{j=1}^k a_jx_j,\qquad
		\bar{y}=\sum_{j=1}^k a_jy_j.
	\end{align*}
	Then
	\begin{equation}\label{v1}
		\begin{aligned}
			f\!\left(N+n-\bar{y}\right)
			&\le f(N)+f(n)-\sum_{j=1}^k a_jf(x_j)
			+\frac{\bar{y}-\bar{x}}{N-n}\bigl(f(n)-f(N)\bigr)\\
			&\quad-c\left[(N+n)(\bar{x}+\bar{y})
			-\left(2Nn+\sum_{j=1}^k a_jx_j^2+\bar{y}^{\,2}\right)\right].
		\end{aligned}
	\end{equation}
\end{thm}

\begin{proof}
	By strong convexity with modulus $c$, for every $t\in[n,N]$,
	\begin{equation}\label{10}
		f(t)\le
		\frac{N-t}{N-n}f(n)+\frac{t-n}{N-n}f(N)
		-c(N-t)(t-n).
	\end{equation}
	Replacing $t$ by $N+n-t$ in \eqref{10} gives
	\begin{equation}\label{11}
		f(N+n-t)\le
		\frac{t-n}{N-n}f(n)+\frac{N-t}{N-n}f(N)
		-c(t-n)(N-t).
	\end{equation}
	Setting $t=x_j$ in \eqref{10}, we obtain
	\begin{equation}\label{12}
		f(x_j)\le
		\frac{N-x_j}{N-n}f(n)+\frac{x_j-n}{N-n}f(N)
		-c\bigl((N+n)x_j-Nn-x_j^2\bigr).
	\end{equation}
	Multiplying \eqref{12} by $a_j$ and summing over
	$j=1,\ldots,k$ yields
	\begin{equation}\label{13}
			\sum_{j=1}^k a_jf(x_j)
			\le \frac{N-\bar{x}}{N-n}f(n)
			+\frac{\bar{x}-n}{N-n}f(N)-c\left[(N+n)\bar{x}
			-Nn-\sum_{j=1}^k a_jx_j^2\right].
	\end{equation}
	Setting $t=\bar{y}$ in \eqref{11} gives
	\begin{equation}\label{14}
	 	f(N+n-\bar{y})
			\le \frac{\bar{y}-n}{N-n}f(n)
			+\frac{N-\bar{y}}{N-n}f(N)-c\left[(N+n)\bar{y}-Nn-\bar{y}^{\,2}\right].
		\end{equation}
	Adding \eqref{13} and \eqref{14} and collecting terms, we find
	\begin{align*}
		\sum_{j=1}^k a_jf(x_j)+f(N+n-\bar{y})
		&\le f(n)+f(N)
		+\frac{\bar{y}-\bar{x}}{N-n}
		\bigl(f(n)-f(N)\bigr)\\
		&\quad-c\left[(N+n)(\bar{x}+\bar{y})
		-2Nn-\sum_{j=1}^k a_jx_j^2-\bar{y}^{\,2}\right].
	\end{align*}
	Rearranging proves \eqref{v1}.
\end{proof}
\begin{remark}	\begin{enumerate}[(1)]
\item[(i)] The correction term in \eqref{v1} is nonnegative because 
\begin{equation}\label{15}	(N+n)\!\left(\bar{x}+\bar{y}\right)	-\left( 2Nn+\sum_{j=1}^k a_j x_j^2+	\left( \bar{y}\right)^2\right) \ge 0.	\end{equation} In support of \eqref{15}, we write
\begin{equation}\label{16} 	\begin{aligned}	(N+n)\!\left(\bar{x} +\bar{y} \right)	-\left( 2Nn+\sum_{j=1}^k a_j x_j^2+
\left(\bar{y} \right) ^2\right)\\	& \hspace{-8cm}=\left[(N+n)\bar{x}-\sum_{j=1}^k a_j x_j^2-nN\right]  +\left[  (N+n)\bar{y}-(\bar{y})^2-nN\right].   	\end{aligned}	\end{equation}
It suffices to verify that the expressions inside the first and second square brackets on the right-hand side of \eqref{16} are nonnegative. Indeed, the first bracket equals $\sum_j a_j(x_j-n)(N-x_j)$, and the second equals $(\bar y-n)(N-\bar y)$. Both are nonnegative.\item[(ii)]
For a convex function, the same proof with $c=0$ gives 	\begin{equation}\label{17}	f\!\left( N+n-\bar{y}\right)	\le f(N)+f(n)-\sum_{j=1}^k a_j f(x_j)	+\frac{\bar{y} -\bar{x}}{N-n}\bigl( f(n)-f(N) \bigr),\end{equation} which recovers the generalized Jensen--Mercer inequality in \cite{mohebbi24}.	\item[(iii)] If we choose $\bar{x}=\bar{y}$, then \eqref{v1} gives 
\begin{equation}\label{v11} 
	\begin{aligned} 
		f(N+n-\bar{x}) 
		&\le f(N)+f(n)-\sum_{j=1}^k a_jf(x_j)\\ 
		&\quad-c\left[ 
		2(N+n)\bar{x} 
		-\left(2Nn+\sum_{j=1}^k a_jx_j^2+\bar{x}^{\,2}\right) 
		\right]. 
	\end{aligned} 
\end{equation} 
\end{enumerate} 

We compare \eqref{v11} with \cite[Theorem~6]{Ivelic2024a}. 
In the notation of that theorem, set 
$$ 
m=n,\quad M=N,\quad A_k=\sum_{j=1}^k a_j=1,\quad 
d=N+n-\bar{x},\quad 
\lambda_j=\frac{N-x_j}{N-n}. 
$$ 
These choices satisfy its hypotheses, since $d,x_j\in[n,N]$ and 
$\lambda_j\in[0,1]$ for $j=1,\ldots,k$. The resulting correction term is 
\[ 
2c(N-n)^2\sum_{j=1}^k a_j\lambda_j(1-\lambda_j) 
+c\sum_{j=1}^k a_j(x_j-\bar{x})^2. 
\] 
Since 
\begin{equation}\label{new2} 
\begin{aligned} 
	&2(N-n)^2\sum_{j=1}^k a_j\lambda_j(1-\lambda_j) 
	+\sum_{j=1}^k a_j(x_j-\bar{x})^2\\ 
	&\quad=2(N-\bar{x})(\bar{x}-n) 
	-\sum_{j=1}^k a_j(x_j-\bar{x})^2\\ 
	&\quad=2(N+n)\bar{x} 
	-\left(2Nn+\sum_{j=1}^k a_jx_j^2+\bar{x}^{\,2}\right), 
\end{aligned} 
\end{equation} 
this specialization of \cite[Theorem~6]{Ivelic2024a} 
coincides with \eqref{v11}, which also coincides with our Theorem~\ref{JMstrong}.\\ 
Dropping the nonnegative term 
$2c(N-n)^2\sum_{j=1}^k a_j\lambda_j(1-\lambda_j)$ 
from its correction term gives the weaker bound 
\begin{equation}\label{new1} 
f(N+n-\bar{x}) 
\le f(N)+f(n)-\sum_{j=1}^k a_jf(x_j) 
-c\sum_{j=1}^k a_j(x_j-\bar{x})^2. 
\end{equation} 
Furthermore, 
\begin{equation}\label{new3} 
\sum_{j=1}^k a_j(x_j-\bar{x})^2 
\le (N-\bar{x})(\bar{x}-n). 
\end{equation} 
Consequently, \eqref{v11} implies 
\begin{equation}\label{new4} 
f(N+n-\bar{x}) 
\le f(N)+f(n)-\sum_{j=1}^k a_jf(x_j) 
-c(N-\bar{x})(\bar{x}-n). 
\end{equation} 
By \eqref{new3}, the correction term in \eqref{new4} is at least as large as 
that in \eqref{new1}. Both \eqref{new1} and \eqref{new4} 
are consequences of the stronger bound \eqref{v11}.
\end{remark}\vspace{0.2cm}
Next, we present the reverse inequality of \eqref{6} for strongly convex functions.\vspace{0.2cm}
\begin{thm}\label{t2}
Let $I$ be an open interval containing $[n,N]$, and let $f:I\to\mathbb{R}$ be differentiable and strongly convex on $[n,N]$ with modulus $c>0$. Let $n<N$ and $n\le t\le N$. Then \begin{equation}\label{e1}	f(N)+f(n)-f(t)	\le f(N+n-t)+\frac{N-n}{2}\bigl(f'(N)-f'(n)\bigr)-c\left( (N-t)^2+(n-t)^2\right).	\end{equation}
\end{thm}
\begin{proof}
Since $g(u)=f(u)-cu^2$ is convex on $[n,N]$, its supporting-line inequality gives
\begin{equation}\label{18}	f(y) \le f(x)+f'(y)(y-x)-c(x-y)^2,	\end{equation}	for any $x,y \in [n,N]$. Substituting $x = t$ and $y = N$ in \eqref{18}, we obtain\begin{equation}\label{19}	f(N)\le f(t)+f'(N)(N-t)-c(t-N)^2.\end{equation}Again, setting $x = t$ and $y = n$ in \eqref{18}, we have\begin{equation}\label{20}	f(n)\le f(t)+f'(n)(n-t)-c(t-n)^2.
\end{equation}Choosing $x = N + n - t$ and $y = N$ in \eqref{18}, we obtain\begin{equation}\label{21}	f(N)\le f(N+n-t)+f'(N)(t-n)-c(n-t)^2.\end{equation}If $x=N+n-t$ and $y=n$ in \eqref{18}, we get	\begin{equation}\label{22}	f(n)\le f(N+n-t)+f'(n)(t-N)-c(N-t)^2.\end{equation}Adding \eqref{19}--\eqref{22} gives $2[f(N)+f(n)]\le 2f(t)+2f(N+n-t)+(N-n)[f'(N)-f'(n)]-2c[(N-t)^2+(n-t)^2]$. Dividing by $2$ and rearranging proves \eqref{e1}.
\end{proof}
As a consequence of Theorem \ref{t2}, we obtain the following corollary, which is a complementary inequality to \eqref{5}. \vspace{0.2cm}
\begin{corollary}\label{cr2} Let $I$ be an open interval containing $[n,N]$ and let $f:I\to\mathbb{R}$ be differentiable and strongly convex on $[n,N]$ with modulus $c>0$. Let $n<N$ and
	$n\le x_j \le N$, and let $a_j$ be positive scalars such that $\sum_{j=1}^k a_j = 1, \bar{x}=\sum_{j=1}^{k}a_jx_j$. Then	\begin{equation}\label{vidyanshu1} \begin{aligned} \hspace{0cm}f(N)+f(n)-\sum_{j=1}^k a_j f(x_j) \le f\!\left(N+n-\bar{x}\right)	+\frac{N-n}{2}\bigl(f'(N)-f'(n)\bigr)\\	&\hspace{-7.5cm}-c\left[ (N-\bar{x})^2+(n-\bar{x})^2+\sum_{j=1}^{k} a_j(x_j-\bar{x})^2\right].\end{aligned}\end{equation}
\begin{proof}
By \eqref{3}, $f(\bar x)\le\sum_j a_jf(x_j)-c\sum_j a_j(x_j-\bar x)^2$. Substitute this estimate for $f(\bar x)$ into \eqref{e1} with $t=\bar x$ and rearrange to obtain the result.
\end{proof}
\end{corollary}
The well-known arithmetic mean--Heinz mean inequality \cite{hardy52} for any $a, b > 0$ states that\begin{equation}\label{neq5} H_{Z_t}(a,b)\le A(a,b), \end{equation} where \begin{equation*} A(a,b)=\frac{a+b}{2}, \qquad	H_{Z_t}(a,b)=\frac{a^{1-t}b^{t}+a^{t}b^{1-t}}{2},		\quad 0\le t\le 1.	\end{equation*}
We now give a counterpart of \eqref{neq5}.\vspace{0.2cm}
\begin{corollary}\label{cr3}
	Let $a$ and $b$ be distinct positive numbers, and set
	$
	c=\frac{1}{2}\left(\ln\frac{b}{a}\right)^2
	\min\left\{1,\frac{b}{a}\right\}.
	$
	Then, for $0\le t\le1$,
	\begin{equation}\label{neq4}
		H_{Z_t}(a,b)
		\ge A(a,b)-\frac14(b-a)(\ln b-\ln a)
		+\frac{ca}{2}\bigl[(1-t)^2+t^2\bigr].
	\end{equation}
\end{corollary}
\begin{proof}
	Set $x=b/a>0$. Since $a\ne b$, we have $x\ne1$. Let
	$f(t)=x^t$ on the open interval $I=\mathbb{R}$. On $[0,1]$,
	$f$ is strongly convex with modulus
	$c=\frac12(\ln x)^2\min\{1,x\}$.
	Applying Theorem \ref{t2} with $n=0$ and $N=1$, we obtain
	\begin{equation}\label{neq2}
		f(1)+f(0)
		\le f(t)+f(1-t)
		+\frac12\bigl(f'(1)-f'(0)\bigr)
		-c\bigl[(1-t)^2+t^2\bigr],
		\quad 0\le t\le1.
	\end{equation}
	Since $f'(t)=(\ln x)x^t$, substituting $f(t)=x^t$ into
	\eqref{neq2} and dividing by $2$ gives
	\begin{equation}\label{neq3}
		\frac{x+1}{2}
		\le \frac{x^t+x^{1-t}}{2}
		+\frac{(x-1)\ln x}{4}
		-\frac{c}{2}\bigl[(1-t)^2+t^2\bigr],
		\quad 0\le t\le1.
	\end{equation}
	Multiplying \eqref{neq3} by $a$ and using $x=b/a$ yields
	\eqref{neq4}.
\end{proof}
It was shown in \cite{kian13} that if $f:I\to\mathbb{R}$ is convex, then for any $n\le x<y\le N$,
\begin{equation*}f\!\left(N+n-\frac{x+y}{2}\right)\le	f(N)+f(n)-\int_{0}^{1} f\bigl((1-t)x+ty\bigr)\,dt\end{equation*}
and	\begin{equation*} \frac{1}{y-x}\int_{x}^{y} f(N+n-t)\,dt\le f(N)+f(n)-\frac{f(x)+f(y)}{2}. \end{equation*}
The following result gives reverse inequalities for strongly convex functions.\vspace{0.2cm}
\begin{thm}\label{t3} Let $I$ be an open interval containing $[n,N]$ and let $f:I\to\mathbb{R}$ be differentiable and strongly convex on $[n,N]$ with modulus $c>0$. Let $n\le x<y\le N$. Then \begin{equation}\label{e3} \begin{aligned}
\hspace{0cm} f\!\left(N+n-\frac{x+y}{2}\right) \ge f(N)+f(n)-\int_{0}^{1}f\bigl((1-t)x+ty\bigr)\,dt-\frac{N-n}{2}\bigl(f'(N)-f'(n)\bigr)\\ & \hspace{-10cm}+c\left[\left( N-\frac{x+y}{2}\right) ^2+\left( n-\frac{x+y}{2}\right) ^2+\frac{(x-y)^2}{12}\right]
\end{aligned}\end{equation}and\begin{equation}\label{e4}	\begin{aligned} \hspace{0cm} \frac{1}{y-x}\int_{x}^{y} f(N+n-t)\,dt \ge f(N)+f(n)-\frac{f(x)+f(y)}{2}-\frac{N-n}{2}\bigl(f'(N)-f'(n)\bigr)\\	&\hspace{-12cm}+\frac{c}{y-x}\left[\frac{(N-x)^3}{3}-\frac{(N-y)^3}{3}+\frac{(n-x)^3}{3}-\frac{(n-y)^3}{3}+\frac{(y-x)^3}{6}\right].
\end{aligned} \end{equation}
\end{thm}
\begin{proof} To prove \eqref{e3}, we first apply the Hermite--Hadamard inequality \eqref{8} on $[x,y]$ and then apply Theorem \ref{t2} with $t=\frac{x+y}{2}$, that is, \begin{equation*} \begin{aligned}
\hspace{0cm}f(N)+f(n)-\int_{0}^{1} f\bigl((1-t)x+ty\bigr)\,dt\le f(N)+f(n)- f\left( \frac{x+y}{2}\right)-\frac{c}{12}(x-y)^2\\	& \hspace{-8.5cm} \le f\!\left(N+n-\frac{x+y}{2}\right)
+\frac{N-n}{2}\bigl(f'(N)-f'(n)\bigr)\\	&\hspace{-8.5cm}-c\left[\left(N-\frac{x+y}{2}\right)^2+\left(n-\frac{x+y}{2}\right)^2+\frac{(x-y)^2}{12}\right]. 	\end{aligned}
\end{equation*}	Applying \eqref{8} on $[x,y]$, we get\begin{equation}\label{e5}
f(N)+f(n)-\frac{f(x)+f(y)}{2} \le f(N)+f(n)-\frac{1}{y-x}\int_{x}^{y} f(t)\,dt-\frac{c}{6}\left( x-y\right)^2.
\end{equation}
Applying Theorem \ref{t2}, for $n \le x \le t \le y \le N $ and integrating with respect to $t$ from $x$ to $y$, we find that
\begin{equation}\label{e6} \begin{aligned}	f(N)+f(n)-\frac{1}{y-x}\int_{x}^{y} f(t)\,dt \le \frac{1}{y-x}\int_{x}^{y} f(N+n-t)\,dt+\frac{N-n}{2}\bigl(f'(N)-f'(n)\bigr)\\	&\hspace{-11cm}-\frac{c}{y-x}\left[\frac{(N-x)^3}{3}-\frac{(N-y)^3}{3}+\frac{(n-x)^3}{3}-\frac{(n-y)^3}{3}\right].	\end{aligned}
\end{equation}	The inequality \eqref{e4} follows by combining \eqref{e5} and \eqref{e6}.
\end{proof}
Next, we present the counterparts of the Jensen--Mercer type inequality \eqref{v1}.\vspace{0.2cm}
\begin{thm}\label{t4} Suppose that the assumptions of Theorem \ref{t1} hold and that $f$ is differentiable on $I$. Then\begin{equation}\label{35} \begin{aligned}
f\!\left(N+n-\bar{y}\right) \ge f(N)+f(n)-\sum_{j=1}^k a_j f(x_j)\\	&\hspace{-5cm}-\frac{N-n}{2}\bigl(f'(N)-f'(n)\bigr)-\left(\frac{\bar{y}-\bar{x}}{2}\right)\bigl(f'(N)+f'(n)\bigr)\\&\hspace{-5cm}+\frac{c}{2}\left[2\sum_{j=1}^{k}a_j x_j^2+(N^2+n^2)-2(n+N)\bar{x}+(n-\bar{y})^2+(N-\bar{y})^2\right]. 	\end{aligned}\end{equation}
\end{thm}
\begin{proof}
Adding \eqref{19} and \eqref{20}, we get\begin{equation}\label{36}
f(N)+f(n)\le 2f(t)+f'(N)(N-t)+f'(n)(n-t)-c\left[ (t-n)^2+(t-N)^2\right].\end{equation}Putting $t = x_j $ $(j = 1, 2, \ldots, k)$ in \eqref{36}, multiplying by $ a_j$ $(j = 1, 2, \ldots, k)$, and then summing the resulting inequalities from $j = 1$ to $k$, we get\begin{equation}\label{37}\begin{aligned}	\hspace{-5cm}f(N)+f(n)
\le	2\sum_{j=1}^k a_j f(x_j)	+f'(N)\!\left(N-\bar{x}\right)	+f'(n)\!\left(n-\bar{x}\right)\\	& \hspace{-9cm}-c\left[2\sum_{j=1}^{k}a_j x_j^2+(N^2+n^2)-2(n+N)\bar{x}\right].\end{aligned}
\end{equation}Replacing $t$ by $N+n-\bar{y}$ in \eqref{36}, we get
\begin{equation}\label{38}
f(N)+f(n)	\le	2f\!\left(N+n-\bar{y}\right)+f'(N)\!\left(\bar{y}-n\right)
+f'(n)\!\left(\bar{y}-N\right)-c\left[(n-\bar{y})^2+(N-\bar{y})^2\right].\end{equation}
Adding \eqref{37} and \eqref{38}, dividing by $2$, and rearranging gives \eqref{35}.\end{proof}
\begin{remark}\begin{enumerate}[(1)]	\item[(i)] For $\bar{x}=\bar{y}$, inequality \eqref{35} agrees with Corollary \ref{cr2}.\vspace{0.1cm}	\item[(ii)] For a differentiable convex function, the same argument with $c=0$ gives\begin{equation}\label{39} \begin{aligned} \hspace{0cm}f(N)+f(n)-\sum_{j=1}^k a_j f(x_j)\le f\!\left(N+n-\bar{y}\right)+\frac{N-n}{2}\bigl(f'(N)-f'(n)\bigr)\\	&\hspace{-7cm} +\left(\frac{\bar{y}-\bar{x}}{2}\right)	\bigl(f'(N)+f'(n)\bigr),\\ \end{aligned} \end{equation}	which agrees with the result in \cite{mohebbi24}. \end{enumerate}
\end{remark}\vspace{0.2cm}
Now, we present a lemma that will be used in the proof of Theorem \ref{t5}.\vspace{0.2cm} 
\begin{lemma}\label{l4}
	Let $f:I\to\mathbb{R}$ be differentiable and strongly convex
	with modulus $c>0$. Let $n,N\in I$ with $n<N$, and assume that
	$f'\in S[n,N]$, the space of Riemann integrable functions on $[n,N]$.
	For $x\in[n,N]$, define
	\begin{equation}\label{neq6}
		\delta_{n,N}(x)
		=\frac{(N-x)^2I_1(x)+(n-x)^2I_2(x)}{2(N-n)},
	\end{equation}
	where
	\begin{align*}
		I_1(x)
		&=\int_0^1 v\bigl[f'(n+v(N-x))
		-f'(N+v(x-N))\bigr]\,dv,\\
		I_2(x)
		&=\int_0^1 v\bigl[f'(n+v(x-n))
		-f'(N+v(n-x))\bigr]\,dv.
	\end{align*}
	Then
	$$
	\delta_{n,N}(x)
	=\frac{f(x)+f(N+n-x)}{2}
	-\frac{1}{N-n}\int_n^N f(t)\,dt.
$$
	In particular, $\delta_{n,N}((n+N)/2)\le0$.
\end{lemma}
\begin{proof}
	Fix $z\in[n,N]$. Integration by parts gives
	\begin{align*}
		\int_n^z (t-n)f'(t)\,dt
		&=(z-n)f(z)-\int_n^z f(t)\,dt,\\
		\int_z^N (N-t)f'(t)\,dt
		&=-(N-z)f(z)+\int_z^N f(t)\,dt.
	\end{align*}
	Subtracting the second identity from the first, we obtain
	\begin{equation}\label{eq:integral-identity}
		f(z)-\frac{1}{N-n}\int_n^N f(t)\,dt
		=\frac{1}{N-n}\left[
		\int_n^z(t-n)f'(t)\,dt
		-\int_z^N(N-t)f'(t)\,dt
		\right].
	\end{equation}
	In the first integral, set $t=n+v(z-n)$; in the second, set
	$t=N+v(z-N)$. Substituting into \eqref{eq:integral-identity} gives
	\begin{equation}\label{eq:parameter-identity}
		\begin{aligned}
			f(z)-\frac{1}{N-n}\int_n^N f(t)\,dt
			&=\frac{(z-n)^2}{N-n}
			\int_0^1 vf'(n+v(z-n))\,dv\\
			&\quad-\frac{(z-N)^2}{N-n}
			\int_0^1 vf'(N+v(z-N))\,dv.
		\end{aligned}
	\end{equation}
	This identity also appears in \cite{cerone04}. Applying
	\eqref{eq:parameter-identity} with $z=x$ and $z=N+n-x$,
	then adding the resulting equations and dividing by $2$, gives
$$
	\frac{f(x)+f(N+n-x)}{2}
	-\frac{1}{N-n}\int_n^N f(t)\,dt
	=\frac{(N-x)^2I_1(x)+(n-x)^2I_2(x)}{2(N-n)}
	=\delta_{n,N}(x).
$$
	Finally, setting $x=(n+N)/2$ and applying the
	Hermite--Hadamard inequality to the convex function $f$, we obtain
	$$
	\delta_{n,N}\!\left(\frac{n+N}{2}\right)
	=f\!\left(\frac{n+N}{2}\right)
	-\frac{1}{N-n}\int_n^N f(t)\,dt
	\le 0.
	$$
\end{proof}
We next obtain a Hermite--Hadamard type inequality using
Lemma~\ref{l4}.\vspace{0.2cm}

\begin{thm}\label{t5}
	Let $f:I\to\mathbb{R}$ be differentiable and strongly convex with
	modulus $c>0$, and let $n,N\in I$ with $n<N$. Assume that
	$f'\in S[n,N]$. For every $x\in[n,N]$, we have
	\begin{equation}\label{rm}
		\begin{aligned}
			f\!\left(\frac{n+N}{2}\right)
			&\le \delta_{n,N}(x)
			+\frac{1}{N-n}\int_n^N f(t)\,dt
			-\frac{c}{4}(N+n-2x)^2\\
			&\le \frac{f(n)+f(N)}{2}
			-c\left[\frac{(N+n-2x)^2}{4}
			+(N-x)(x-n)\right],
		\end{aligned}
	\end{equation}
	where $\delta_{n,N}(x)$ is defined in \eqref{neq6}.
\end{thm}

\begin{proof} 
	Applying strong convexity to $x$ and $N+n-x$ with equal weights 
	gives 
	\begin{equation}\label{r1} 
		f\left(\frac{n+N}{2}\right) 
		\le \frac{f(x)+f(N+n-x)}{2} 
		-\frac{c}{4}(N+n-2x)^2. 
	\end{equation} 
	By Lemma~\ref{l4}, we have
	\begin{equation}\label{r2} 
		\frac{f(x)+f(N+n-x)}{2} 
		=\delta_{n,N}(x)+\frac{1}{N-n}\int_n^N f(t)\,dt. 
	\end{equation} 
	The first inequality in \eqref{rm} follows by combining \eqref{r1} and \eqref{r2}, while the second inequality follows from \eqref{7*}. 
\end{proof}
\section{Applications to Strongly $f$-Divergences}\label{Ap}
Nonnegative measures of dissimilarity between probability distributions
are known as divergence measures. An important class of such 
measures is generated by convex functions and is referred to as 
Csisz\'ar $f$-divergences; see \cite{csiszar67} for the original formulation. In recent years, various bounds for several divergence measures have been established; see \cite{dragomir14, sason16}. Furthermore, Iveli\'c \cite{ivelic19} introduced the notion of strongly $f$-divergence, which extends the classical framework to the setting of strongly convex functions.\\ Fix $0<n<N$ and let $f:(0,\infty)\to\mathbb{R}$ be a convex function that is strongly convex with modulus $c>0$ on $[n,N]$. Let 	$a=(a_1,\dots,a_k)$ and $b=(b_1,\dots,b_k)$ 
be probability distributions satisfying $a_j>0,\quad b_j>0, \quad \sum_{j=1}^k a_j=\sum_{j=1}^k b_j=1.$
We define the strongly Csisz\'ar $f$-divergence by	$$	\tilde{D}_f(b,a)	=\sum_{j=1}^{k} a_j 	f\!\left(\frac{b_j}{a_j}\right),$$
and the strongly chi-square divergence by	$$\tilde{D}_{\chi^2}(b,a)	=	\sum_{j=1}^{k} a_j \left(\frac{b_j}{a_j}\right)^2 -1.$$
Csisz\'ar $f$-divergences play a fundamental role in probability theory, statistics, information theory, physics, and computer science; 
see \cite{crooks, dragomir00a, dragomir00b, polyanskiy25}. They unify and generalize many classical divergence measures. Among the most widely used special cases are the Kullback--Leibler divergence, 
$\chi^{2}$-divergence, Hellinger divergence, Bhattacharyya-type divergence, 
and Jeffreys divergence, which will be discussed in detail later in this section	\\
We now apply the results established in the previous section 
to derive upper and lower bounds for the strongly $f$-divergence 	$\tilde{D}_f(b,a)$.\vspace{0.2cm}
\begin{thm}\label{t31} Let $0<n<N$ and assume $\frac{b_j}{a_j}\in[n,N], 
\quad j=1,\dots,k.$
Let $I\subset(0,\infty)$ be an open interval containing $[n,N]$, and let $f:I\to\mathbb R$ be twice continuously differentiable with $f(1)=0$ and modulus on $[n,N]$ given by$$	c=\frac12\inf_{x\in[n,N]} f''(x)>0.$$
Then
\begin{equation}\label{vidyanshu123}	\begin{aligned}& \Delta_f(n,N)-\frac{N-n}{2}(f'(N)-f'(n))+c[(N-1)^2+(n-1)^2]  \\& \le \tilde{D}_f(b,a)	- c\,\tilde{D}_{\chi^2}(b,a)  \\& \le \Delta_f(n,N)-2c(N+n-Nn-1),\end{aligned}
\end{equation}where $\Delta_f(n,N)=f(N)+f(n)-f(n+N-1).$
\end{thm}
\begin{proof} Let $x_j=\frac{b_j}{a_j},\quad j=1,2,\dots,k.$ Since $a$ and $b$ are probability distributions, $\sum_{j=1}^k a_j=\sum_{j=1}^k b_j=1.$ Hence,\begin{equation*}\label{meanvalue}	\bar{x}	=\sum_{j=1}^{k} a_j x_j	=\sum_{j=1}^{k} a_j \frac{b_j}{a_j}	=\sum_{j=1}^{k} b_j	=1.\end{equation*} 	
The function $g(x)=f(x)-cx^2$ is convex on $[n,N]$. Jensen's inequality gives $\sum_j a_jg(x_j)\ge g(1)$. Applying the argument used in the proof of Theorem~\ref{t2} to the convex function $g$ (corresponding to the case $c=0$) at $t=1$ gives \begin{equation}\label{lowerstep}	\Delta_f(n,N)	-\frac{N-n}{2}\big(f'(N)-f'(n)\big)	+c\big[(N-1)^2+(n-1)^2-1\big]\le \sum_{j=1}^{k} a_j f(x_j)-c\sum_{j=1}^{k} a_j x_j^2,\end{equation}
Applying the classical Jensen--Mercer inequality to $g$ gives
\begin{equation}\label{upperstep}
\sum_{j=1}^{k} a_j f(x_j)
-c\sum_{j=1}^{k} a_j x_j^2	\le \Delta_f(n,N)
-c[2(N+n-Nn)-1].\end{equation} 	
Observing that$$\sum_{j=1}^{k} a_j f(x_j)=\tilde{D}_f(b,a),$$	and
$$\sum_{j=1}^{k} a_j x_j^2	=\sum_{j=1}^{k} a_j	\left(\frac{b_j}{a_j}\right)^2=\tilde{D}_{\chi^2}(b,a)+1,$$and using $\bar{x}=1$, inequalities \eqref{lowerstep} and \eqref{upperstep} yield \eqref{vidyanshu123}. This completes the proof.
\end{proof}
We now illustrate Theorem~\ref{t31} for several important divergence measures. In each case, the modulus of strong convexity is computed explicitly.\vspace{0.2cm}
\begin{example}\label{eg1}
	The strongly Kullback--Leibler divergence is defined by
	\begin{equation*}
		\widetilde{KL}(b,a)=\sum_{j=1}^{k} b_j \ln \frac{b_j}{a_j}.
	\end{equation*}
	Let $f(x)=x\ln x$. Then
	$
	f''(x)=\frac{1}{x},\quad x>0.
	$
	Since $\frac{b_j}{a_j}\in[n,N]\subset(0,\infty)$, we have
	$$
	\inf_{x\in[n,N]}f''(x)=\frac{1}{N}.
	$$
	Hence, $f$ is strongly convex on $[n,N]$ with modulus
	$
	c=\frac{1}{2N}.
	$
	Applying Theorem~\ref{t31}, we obtain
	\begin{equation*}
		\begin{aligned}
			& \Delta_f(n,N)-\frac{N-n}{2}\ln\frac{N}{n}
			+\frac{1}{2N}\big[(N-1)^2+(n-1)^2\big] \\
			& \le \widetilde{KL}(b,a)-\frac{1}{2N}\widetilde{D}_{\chi^2}(b,a) \\
			& \le \Delta_f(n,N)-\frac{1}{N}(N+n-Nn-1).
		\end{aligned}
	\end{equation*}
\end{example}
\begin{example}\label{eg2}	The strongly squared Hellinger divergence is defined by
\begin{equation*}	\widetilde{D}_h^2(b,a)=\sum_{j=1}^{k}(\sqrt{a_j}-\sqrt{b_j})^2.
\end{equation*}
Let $f(x)=(\sqrt{x}-1)^2$. Then
$f''(x)=\frac{1}{2x^{3/2}}, \quad x>0.$ Since $\frac{b_j}{a_j}\in[n,N]$, we obtain	$$\inf_{x\in[n,N]} f''(x)=\frac{1}{2N^{3/2}}.$$Hence, $f$ is strongly convex on $[n,N]$ with modulus
$c=\frac{1}{4N^{3/2}}.$
Applying Theorem~\ref{t31}, we obtain	\begin{equation*}	\begin{aligned}	& \Delta_f(n,N)-\frac{N-n}{2}\left(\frac{1}{\sqrt n}-\frac{1}{\sqrt N}\right)	+\frac{1}{4N^{3/2}}\big[(N-1)^2+(n-1)^2\big] \\	& \le \widetilde{D}_h^2(b,a)-\frac{1}{4N^{3/2}}\widetilde{D}_{\chi^2}(b,a) \\	& \le \Delta_f(n,N)-\frac{1}{2N^{3/2}}(N+n-Nn-1).\end{aligned}	\end{equation*}\end{example}
\begin{example}\label{eg3}	The Bhattacharyya affinity deficit, equal to one half of the squared Hellinger divergence, is defined by\begin{equation*}\widetilde{D}_{B}(b,a)=1-\sum_{j=1}^{k}\sqrt{a_j b_j}.	\end{equation*}Let $f(x)=1-\sqrt{x}$. Then $f''(x)=\frac{1}{4x^{3/2}}, \quad x>0.$ Since $\frac{b_j}{a_j}\in[n,N]$, we have$$	\inf_{x\in[n,N]} f''(x)=\frac{1}{4N^{3/2}}.$$Hence, $f$ is strongly convex on $[n,N]$ with modulus $c=\frac{1}{8N^{3/2}}.$ Applying Theorem \ref{t31}, we obtain\begin{equation*}	\begin{aligned}
& \Delta_f(n,N)-\frac{N-n}{4}
\left(\frac{1}{\sqrt{n}}-\frac{1}{\sqrt{N}}\right)+\frac{1}{8N^{3/2}}\big[(N-1)^2+(n-1)^2\big] \\	& \le \widetilde{D}_{B}(b,a)-\frac{1}{8N^{3/2}}\widetilde{D}_{\chi^2}(b,a) \\	& \le \Delta_f(n,N)-\frac{1}{4N^{3/2}}(N+n-Nn-1).	\end{aligned}	\end{equation*}\end{example}
\begin{example}\label{eg4}The strongly Jeffreys divergence is defined by
\begin{equation*}\widetilde{D}_J(b,a)=\sum_{j=1}^{k}(b_j-a_j)\ln\frac{b_j}{a_j}.\end{equation*}Let $f(x)=(x-1)\ln x$. Then $f''(x)=\frac{x+1}{x^2}, \quad x>0.$ Since $\frac{b_j}{a_j}\in[n,N]$, we obtain$$\inf_{x\in[n,N]} f''(x)=\frac{N+1}{N^2}.$$Hence, $f$ is strongly convex on $[n,N]$ with modulus $c=\frac{N+1}{2N^2}.$ Applying Theorem~\ref{t31}, we obtain\begin{equation*}\begin{aligned}	& \Delta_f(n,N)-\frac{N-n}{2}\left(\frac{N-1}{N}-\frac{n-1}{n}	+\ln\frac{N}{n}\right)	+\frac{N+1}{2N^2}\big[(N-1)^2+(n-1)^2\big] \\	& \le \widetilde{D}_J(b,a)	-\frac{N+1}{2N^2}\widetilde{D}_{\chi^2}(b,a) \\	& \le \Delta_f(n,N)	-\frac{N+1}{N^2}(N+n-Nn-1).	\end{aligned}
\end{equation*}
\end{example}
\section*{Conclusion}
In this paper, we established several Jensen--Mercer-type inequalities for strongly convex functions and applied them to strong $f$-divergences. The obtained inequalities contain explicit quadratic correction terms involving the modulus of strong convexity, leading to explicit bounds for several classical divergence measures. These results may also serve as a basis for further investigations in matrix analysis, operator inequalities, and information theory.
\section*{Acknowledgments}

The authors acknowledge the financial support provided by the National Board of Higher Mathematics (NBHM), Department of Atomic Energy (DAE), India (Grant No.~0211/30/2025/NBHM(R.P.)/R\&D-II/9676).
\section*{Declarations}
\subsection*{Funding}
The authors received no financial support for the preparation of this paper apart from the grant acknowledged above.
\subsection*{Conflict of interest} The authors declare that they have no conflict of interest.
\subsection*{Data availability} No data associated with this publication.


\begin{thebibliography}{1}
	\bibitem{butt22}
	Butt, S.~I., Agarwal, P., Yousaf, S., Guirao, J.~L.:
	Generalized fractal Jensen and Jensen--Mercer inequalities for harmonic convex functions with applications.
	\textit{J. Inequal. Appl.} \textbf{2022}, Article 1 (2022)
	
	\bibitem{csiszar67}
	Csisz\'ar, I.:
	On information-type measures of difference of probability distributions and indirect observations.
	\textit{Studia Sci. Math. Hungar.} \textbf{2}, 299--318 (1967)
	
	\bibitem{cerone04}
	Cerone, P., Dragomir, S.~S.:
	Ostrowski-type inequalities for functions whose derivatives satisfy certain convexity assumptions.
	\textit{Demonstr. Math.} \textbf{37}, 299--308 (2004)
	
	\bibitem{cerone10}
	Cerone, P., Dragomir, S.~S.:
	\textit{Mathematical Inequalities: A Perspective.}
	CRC Press, Boca Raton (2010)
	
	\bibitem{crooks}
	Crooks, G.~E.:
	On measures of entropy and information.
	Technical Note \textbf{7}, 1--20 (2017)
	
	\bibitem{dragomir00a}
	Dragomir, S.~S.:
	Inequalities for Csisz\'ar $f$-divergence in information theory.
	Victoria University, Melbourne (2000)
	
	\bibitem{dragomir00b}
	Dragomir, S.~S.:
	Other inequalities for Csisz\'ar divergence and applications.
	Preprint, RGMIA Res. Rep. Coll. (2000)
	
	\bibitem{dragomir14}
	Dragomir, S.~S.:
	A new quantum $f$-divergence for trace class operators in Hilbert spaces.
	\textit{Entropy} \textbf{16}, 5853--5875 (2014)
	
	\bibitem{hardy52}
	Hardy, G.~H., Littlewood, J.~E., P\'olya, G.:
	\textit{Inequalities.}
	Cambridge University Press, Cambridge (1952)
	
	\bibitem{horvath21}
	Horv\'ath, L.:
	Some notes on Jensen--Mercer type inequalities: Extensions and refinements with applications.
	\textit{Math. Inequal. Appl.} \textbf{24}, 1093--1111 (2021)
	
	\bibitem{ivelic11}
	Iveli\'c, S., Matkovi\'c, A., Pe\v{c}ari\'c, J.~E.:
	On a Jensen--Mercer operator inequality.
	\textit{Banach J. Math. Anal.} \textbf{5}, 19--28 (2011)
	
	\bibitem{ivelic19}
	Iveli\'c Bradanovi\'c, S.:
	Sherman's inequality and its converse for strongly convex functions with applications to generalized $f$-divergences.
	\textit{Turk. J. Math.} \textbf{43}, 2680--2696 (2019)
	
	\bibitem{ivelic24}
	Iveli\'c Bradanovi\'c, S.:
	Improvements of Jensen's inequality and its converse for strongly convex functions with applications to strongly $f$-divergences.
	\textit{J. Math. Anal. Appl.} \textbf{531}, Article 127699 (2024)
	
	\bibitem{Ivelic2024a}
	Iveli\'c Bradanovi\'c, S., Lovri\v{c}evi\'c, N.:
	Generalized Jensen and Jensen--Mercer inequalities for strongly convex functions with applications.
	\textit{J. Inequal. Appl.} \textbf{2024}, Article 112 (2024)
	
	\bibitem{jarad23}
	Jarad, F., Sahoo, S.~K., Nisar, K.~S., Trean\c{t}\u{a}, S., Emadifar, H., Botmart, T.:
	New stochastic fractional integrals and related inequalities of Jensen--Mercer and Hermite--Hadamard--Mercer type for convex stochastic processes.
	\textit{J. Inequal. Appl.} \textbf{2023}, Article 51 (2023)
	
	\bibitem{krnic12}
	Krni\'c, M., Lovri\v{c}evi\'c, N., Pe\v{c}ari\'c, J.:
	On some properties of Jensen--Mercer's functional.
	\textit{J. Math. Inequal.} \textbf{6}, 125--139 (2012)
	
	\bibitem{kian13}
	Kian, M., Moslehian, M.~S.:
	Refinements of the operator Jensen--Mercer inequality.
	\textit{Electron. J. Linear Algebra} \textbf{26}, 742--753 (2013)
	
	\bibitem{khan20}
	Khan, A., Husain, Z., Chu, Y.-M.:
	New estimates for Csisz\'ar divergence and Zipf--Mandelbrot entropy via Jensen--Mercer's inequality.
	\textit{Complexity} \textbf{2020}, Article ID 8856502 (2020)
	
	\bibitem{khan22}
	Khan, A.~R., Rubab, F.:
	Mercer-type variants of the Jensen--Steffensen inequality.
	\textit{Rocky Mt. J. Math.} \textbf{52}, 1693--1712 (2022)
	
	\bibitem{mercer03}
	Mercer, A.~M.:
	A variant of Jensen's inequality.
	\textit{J. Inequal. Pure Appl. Math.} \textbf{4}, Article 73 (2003)
	
	\bibitem{merentes10}
	Merentes, N., Nikodem, K.:
	Remarks on strongly convex functions.
	\textit{Aequationes Math.} \textbf{80}, 193--199 (2010)
	
	\bibitem{moradi18}
	Moradi, H.~R., Omidvar, M.~E., Khan, M.~A., Nikodem, K.:
	Around Jensen's inequality for strongly convex functions.
	\textit{Aequationes Math.} \textbf{92}, 25--37 (2018)
	
	\bibitem{mohebbi24}
	Mohebbi, F.~P., Hassani, M., Omidvar, M.~E., Moradi, H.~R., Furuichi, S.:
	Further Jensen--Mercer type inequalities for convex functions.
	\textit{J. Math. Inequal.} \textbf{18}, 719--737 (2024)
	
	\bibitem{niculescu06}
	Niculescu, C.~P., Persson, L.-E.:
	\textit{Convex Functions and Their Applications.}
	Springer, New York (2006)
	
	\bibitem{nikodem14}
	Nikodem, K.:
	On strongly convex functions and related classes of functions.
	In: Rassias, T.~M. (ed.) \textit{Handbook of Functional Equations: Functional Inequalities}, Vol.~95.
	Springer, New York (2014)
	
	\bibitem{polyak66}
	Polyak, B.~T.:
	Existence theorems and convergence of minimizing sequences for extremal problems with constraints.
	\textit{Dokl. Akad. Nauk SSSR} \textbf{166}, 287--290 (1966)
	
	\bibitem{pecaric92}
	Pe\v{c}ari\'c, J.~E., Proschan, F., Tong, Y.~L.:
	\textit{Convex Functions, Partial Orderings, and Statistical Applications.}
	Academic Press, Boston (1992)
	
	\bibitem{polyanskiy25}
	Polyanskiy, Y., Wu, Y.:
	\textit{Information Theory: From Coding to Learning.}
	Cambridge University Press, Cambridge (2025)
	
	\bibitem{roberts73}
	Roberts, A.~W., Varberg, D.~E.:
	\textit{Convex Functions.}
	Academic Press, New York (1973)
	
	\bibitem{sason16}
	Sason, I., Verd\'u, S.:
	$f$-divergence inequalities.
	\textit{IEEE Trans. Inf. Theory} \textbf{62}, 5973--6006 (2016)
	
	\bibitem{sayyari23}
	Sayyari, Y., Barsam, H.:
	Jensen--Mercer inequality for uniformly convex functions with some applications.
	\textit{Afr. Mat.} \textbf{34}, Article 38 (2023)
	\end{thebibliography}
\end{document}